\documentclass[12pt,a5paper]{article}

\usepackage{amsmath}  
\allowdisplaybreaks
\usepackage{defns,ajr}
\newcommand{\pbc}{\textsc{p}t\textsc{bc}}

\newcommand{\E}[2]{$#1\times10^{#2}$}

\title{Higher order accuracy in the gap-tooth scheme for large-scale 
solutions using microscopic simulators}

\author{A.~J.~Roberts\thanks{Dept.  Maths \& Computing, University of
Southern Queensland, Toowoomba, \textsc{Australia}.
\protect\url{mailto:aroberts@usq.edu.au}} \and
I.~G.~Kevrekidis\thanks{Program in Applied and Computational
Mathematics, Princeton University, Princeton, NJ~08544, USA.
\protect\url{mailto:yannis@Princeton.edu}}}

\begin{document}

\maketitle

\begin{abstract}
We are developing a framework for multiscale computation which enables
models at a ``microscopic'' level of description, for example Lattice
Boltzmann, Monte Carlo or Molecular Dynamics simulators, to perform
modelling tasks at the ``macroscopic'' length scales of interest.
The plan is to use the microscopic rules restricted to small patches of
the domain, the ``teeth'', followed by interpolation to estimate
macroscopic fields in the ``gaps''.
The challenge addressed here is to find general boundary conditions for
the patches of microscopic simulators that appropriately connect the
widely separated ``teeth'' to achieve high order accuracy over the
macroscale.
Here we start exploring the issues in the simplest case when the
microscopic simulator is the quintessential example of a partial
differential equation.
For this case analytic solutions provide comparisons.
We argue that classic high-order interpolation provides patch boundary
conditions which achieve arbitrarily high-order consistency in the
gap-tooth scheme, and with care are numerically stable.
The high-order consistency is demonstrated on a class of linear partial
differential equations in two ways: firstly using the dynamical systems
approach of holistic discretisation; and secondly through the
eigenvalues of selected numerical problems.
When applied to patches of microscopic simulations these patch boundary
conditions should achieve efficient macroscale simulation.
\end{abstract}

\tableofcontents

\section{Introduction}
\label{sec:intro}

\begin{figure}
    \centering
        \includegraphics[width=0.85\textwidth]{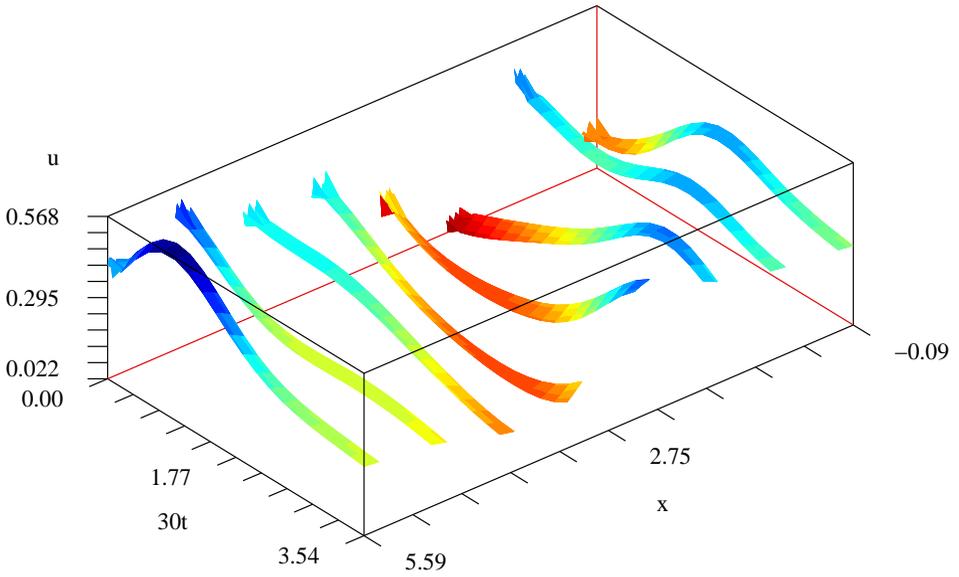} 
	\caption{Gap-tooth solution of Burgers' equation~(\ref{eq:burg})
	on~$[0,2\pi]$ through microsimulation on 8~teeth, each of small
	width~$\pi/20$; the teeth are coupled by special patch boundary
	conditions.  }
    \label{fig:burg3}
\end{figure}

As a preliminary illustration of the gap-tooth scheme~\cite{Gear03,
Samaey03a, Samaey03b}, consider simulating the diffusion and nonlinear
advection of the viscous Burgers' equation
\begin{equation}
    \D tu+100\,u\D xu=\DD xu\,.
    \label{eq:burg}
\end{equation}
Suppose our aim is to simulate the evolution of fields~$u(x,t)$
periodic in~$x$ on the macroscopic length scale~$2\pi$.
See in Figure~\ref{fig:burg3} the continuous time evolution on $m=8$
grid ``points'' in space with macroscopic spacing $H=\pi/4$\,.
However, each ``point'' is actually a microscopic patch of width
$h=\pi/20$\,.
Further, the \emph{only} knowledge that the macroscopic evolution has
of Burgers' \pde~(\ref{eq:burg}) is through the detailed simulation
of the \pde{} within each patch; here we obtain this local detailed
simulation via a discretisation of~(\ref{eq:burg}) on a microscopic
spatial grid of $n=11$ points within each patch, $\Delta x=0.0175$\,,
and on a microscopic time step of $\Delta t\approx 10^{-4}$\,.
This fine scale discretisation of Burgers' \pde~(\ref{eq:burg})
represents a finely detailed model or particle simulation that is too
expensive to use over the entire macroscopic domain.
Our task here is to begin to show how such microscopic simulations in
\emph{relatively small} patches of space may be coupled by patch
boundary conditions derived in Section~\ref{sec:2} to ensure high order
accuracy over the macroscopic~domain.

\begin{figure}
    \centering
    \begin{tabular}{cc}
        \raisebox{20ex}{$u$} & 
        \includegraphics[width=0.85\textwidth]{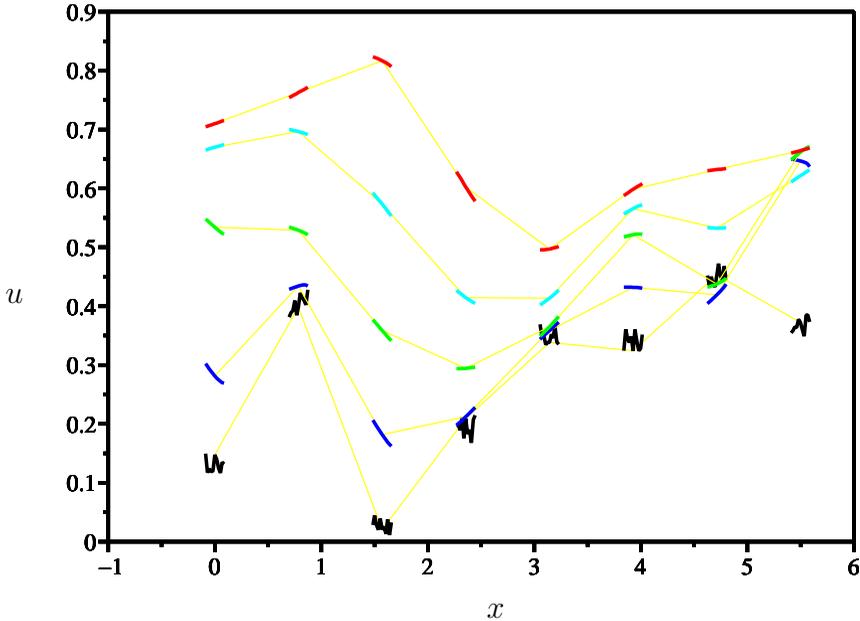}  \\
        & $x$
    \end{tabular}
	\caption{Gap-tooth solution of Burgers' equation~(\ref{eq:burg})
	on~$[0,2\pi]$ on 8~teeth each of small width~$\pi/20$ and coupled
	by special patch boundary conditions.  Solutions $u(x,t)+4t$ are
	plotted at five times $t=0:0.025:0.1$ in different colours and
	connected by yellow lines and with the vertical displacement
	of~$4t$ to help distinguish the plots.}
    \label{fig:burg}
\end{figure}
The example of Figures~\ref{fig:burg3} and~\ref{fig:burg} shows us 
that there are two time scales in the simulation.
Rapidly, the initial internal structure within each tooth (black curves
in Figure~\ref{fig:burg}) smooths by diffusion on the microscopic
time-scale to a local quasi-equilibrium (blue curves).
Then, over longer times the inter-patch coupling exchanges
information between the teeth to guide how the local quasi-equilibria
evolve over macroscopic times.
See the dynamics in Figure~\ref{fig:burg}: the broad hump initially
centred around $x\approx4$ is nonlinearly advected to the right to wrap
around to about $x\approx1$ at the end (red curves); whereas the short
hump initially at $x\approx1$ is dissipated quickly against the
slow moving region at $x\approx2$\,.
Because of the two time scales, we plan future research to implement
``coarse grained'' integration~\cite{Cisternas03, Gear02} which uses
just short bursts of microscopic integration to then extrapolate over a
macroscopic time step.
The result will then be a scheme where the microscopic simulations are
only needed for relatively small patches in space-time.
However, here we concentrate on only the issue of the macroscopic
coupling of small patches across space.

The method of ``holistic discretisation'', developed by Roberts and
Mackenzie~\cite{Roberts98a, Roberts00a, Roberts03c, Mackenzie03},
creates discretisations on a macroscopic grid using systematically
obtained analytic approximations for the subgrid field.
The analytic solutions of Section~\ref{sec:3} in this method are
analogous to the microscopic system simulators in the gap-tooth scheme:
they both provide microscopic solutions which are macroscopically
coupled to neighbouring elements.
This dynamical systems approach adapts the patch boundary conditions to
support the modelling by centre manifold theory~\cite{Roberts00a}.
Then the equivalent \pde\ of the macroscopic dynamic model is found, in
Section~\ref{sec:3}, to confirm the high order consistency for a wide
class of linear \pde{}s.

Lastly, in Section~\ref{sec:4}, we consider a numerical time integrator
for the diffusion equation on patches.  The eigenvalues of the
integrator again confirm the high order accuracy of the proposed patch
boundary conditions.

\section{Couple the patches}
\label{sec:2}
In this section we develop a coupling of the internal dynamics of
patches with their neighbours to achieve high order consistency.
We construct a boundary condition for the flux on the edge of the
microscopic patches that is a natural interpolation of the surrounding
macroscopic field.
The number of required boundary conditions will depend upon the
microscopic simulator~\cite[e.g.]{MacKenzie00a}, and possibilities
other than the flux remain to be explored.
But here we know that boundary conditions on the flux should give rise to
well-posed diffusion-like dynamics.

We introduce the notation in which we typically use capital letters for
macroscopic quantities and lower case letters for microscopic
quantities.
Thus let each of $m$~patches be centred on equi-spaced grid points
$x=X_j=jH$ seen in Figures~\ref{fig:burg3} and~\ref{fig:burg}.
Let each patch be of width~$h$.
Then the edge of a patch is a distance~$h/2$ from its grid point, a
fraction $r=h/(2H)$ to the neighbouring grid point: when $r=\half$ the
neighbouring patches meet and there would be no gap, as in holistic
discretisation~\cite{Roberts98a}.
Here we expect the fraction~$r$ to be small so that the patches are a
relatively small part of the physical domain.
For example, $r=1/10$ in Figures~\ref{fig:burg3} and~\ref{fig:burg}.
Now let $v_j(x,t)$~be the microscopic field in the $j$th~patch.

We use the following identities for discrete operators~\cite{npl61} on
a step size of the macroscopic grid and are careful whether we are
using as a step of~$H$ in~$x$ or a step of~$1$ in~$j$.
In terms of the shift operator, $Ev(x,t)=v(x+H,t)$ or equivalently
$EU_j=U_{j+1}$:
\begin{eqnarray}
    \mbox{centred mean}&& \mu=\half(E^{1/2}+E^{-1/2})\,,
    \label{eq:mean}\\
    \mbox{centred difference}&& \delta=E^{1/2}-E^{-1/2}\,,
    \label{eq:diff}\\
    \mbox{translate/shift}&& E=1+\mu\delta+\half\delta^2\,,
    \label{eq:tran}\\
    \mbox{derivative in~$x$}&& H\partial_x =2\sinh^{-1}\half\delta
    \,,
    \label{eq:deriv}\\
    \mbox{an identity}&& \mu^2=1+\rat14\delta^2\,.
    \label{eq:iden}
\end{eqnarray}
For example, the derivative of the microscopic field on the edge of a
patch, $H\D x{v_j}$ at~$(X_j\pm rH,t)$, may be obtained from~$v_j$
through applying the operator
\begin{eqnarray}
    E^{\pm r}H\partial_x 
    &=&(1+\mu\delta+\half\delta^2)^{\pm r}2\sinh^{-1}\half\delta
    \quad\text{by~(\ref{eq:tran}) and~(\ref{eq:deriv})}
    \nonumber\\
    &=&[1\pm r\mu\delta +\Ord{\delta^2}][\delta+\Ord{\delta^3}]
    \nonumber\\
    &=&\delta\pm r\mu\delta^2 +\Ord{\delta^3}
    \nonumber\\
    &=&\mu\delta \pm r\delta^2 +\Ord{\delta^3}
    \quad\text{by~(\ref{eq:iden}).}
    \label{eq:aoper}
\end{eqnarray}
This last operator just involves evaluation at the grid points~$X_j$
and hence is evaluated from the macroscopic grid values~$U_j$.
This provides the same approximation for the microscopic gradient as
obtained by quadratic interpolation through the neighbouring
macroscopic grid values~\cite[e.g.]{Gear03}.
We proceed to modify such a patch boundary condition in order to obtain
higher order consistency with the surrounding macroscopic variations.

For arbitrary order consistency, as the macroscopic grid size $H\to0$
or as the gradients become small, repeat the previous analysis but
retain more terms, and using~(\ref{eq:iden}) to replace~$\mu^2$ terms:
\begin{eqnarray}
    E^{\pm r}H\partial_x 
    &=&(1+\mu\delta+\half\delta^2)^{\pm r}2\sinh^{-1}\half\delta
    \nonumber\\
    &=&\frac\mu{\sqrt{1+\rat14\delta^2}}
    (1+\mu\delta+\half\delta^2)^{\pm r}2\sinh^{-1}\half\delta
    \nonumber\\
    &=&\mu\delta \pm r\delta^2 
    -(\rat16-\rat12r^2)\mu\delta^3
    \mp r(\rat1{12}-\rat16r^2)\delta^4
    \nonumber\\&&{}
    +(\rat1{30}-\rat18r^2+\rat1{24}r^4)\mu\delta^5
    \pm r(\rat1{90}-\rat1{36}r^2+\rat1{120}r^4)\delta^6
    \nonumber\\&&{}
    -(\rat1{140} -\rat7{240}r^2 +\rat1{72}r^4 -\rat1{720}r^6 )\mu\delta^7
    \nonumber\\&&{}
    \mp r(\rat1{560} -\rat7{1440}r^2 +\rat1{480}r^4 -\rat1{5040}r^6 )\delta^8
    +\Ord{\delta^9}.
    \label{eq:boper}
\end{eqnarray}
Numerical eigenanalysis of the diffusion equation~(\ref{eq:diffn})
reported in Section~\ref{sec:4} confirms the high order accuracy and
stability of the resultant integration scheme with patch boundary
conditions from the above operator.

\section{Achieve high order consistency}
\label{sec:3}

Here we demonstrate analytically that appropriate patch boundary
conditions achieve high order consistency for a wide class of \pde{}s.
Consider the linear \pde
\begin{equation}
    \D tu=\DD xu -c\D xu -b\DDD xu -a\DDDD xu\,,
    \label{eq:lpde}
\end{equation}
for some constants $a$, $b$~and~$c$---we have chosen time and space
scales so that the coefficient of the diffusion term is~$1$.
Specially crafted boundary conditions on small patches ensures
macroscopic consistency.

Following the dynamical systems approach of holistic
discretisation~\cite{Roberts98a, Roberts00a} we introduce the
parameter~$\gamma$ to control the coupling between patches: when
$\gamma=0$ the patches are uncoupled to provide a base for us to apply
centre manifold theory; but when we subsequently set $\gamma=1$ we
recover a dynamical model for the original \pde.
Modify the operator~(\ref{eq:boper}) to invoke the patch boundary
condition (\pbc) that on $x=X_j$ (noting that the $E^{\pm r}$ implies
the left-hand side is evaluated on the edge of the patch at $x=X_j\pm
rH$):
\begin{eqnarray}\hspace*{-1em}
    E^{\pm r}H\partial_x v_j
    &=& \left\{\gamma\left[ \mu\delta \pm r\delta^2 \right]
    +\gamma^2\left[
    -(\rat16-\rat12r^2)\mu\delta^3
    \mp r(\rat1{12}-\rat16r^2)\delta^4 \right]
    \right.\nonumber\\&&{}
    +\gamma^3\left[
    +(\rat1{30}-\rat18r^2+\rat1{24}r^4)\mu\delta^5
    \pm r(\rat1{90}-\rat1{36}r^2+\rat1{120}r^4)\delta^6 \right]
    \nonumber\\&&{}
    +\gamma^4\left[
    -(\rat1{140} -\rat7{240}r^2 +\rat1{72}r^4 -\rat1{720}r^6 )\mu\delta^7
    \right.\nonumber\\&&\left.\left.\quad{}
    \mp r(\rat1{560} -\rat7{1440}r^2 +\rat1{480}r^4 -\rat1{5040}r^6 )\delta^8
   \right] \right\} U_j\,.
   \label{eq:pbcs}
\end{eqnarray}
these \pbc{}s that when $\gamma=0$ the small patches are decoupled and
the resulting insulating boundary conditions, $E^{\pm r}\partial_x
v_j=0$\,, cause the dissipative dynamics of~(\ref{eq:lpde}) in each
patch to decay exponentially quickly to some constant field in each
patch, namely $v_j(x,t)\to U_j$ for each of the
$m$~patches.\footnote{Why are we only using one pair of boundary
conditions for the apparently fourth order \pde~(\ref{eq:lpde})? One
answer is that the \pde\ may be viewed as an equivalent \pde\ for a
microscopic simulator that only requires one pair of boundary
conditions.
For example, if the microscopic simulator is simply a fine scale
discretisation of a \pde, such as $u_t =-c(\delta/h)u
+(\delta^2/h^2)u$\,, then only one pair of boundary conditions are
needed for the simulator, but it has a high order equivalent \pde\ such
as~(\ref{eq:lpde}).
Another answer is that there is no physical boundary at the edge of a
patch and so we only need resolve smooth subgrid fields and for smooth
solutions we need only treat the higher order terms as perturbations;
see that our error terms are expressed as $\Ord{a^2+b^2+c^2}$ for this
reason.} For non-zero coupling parameter~$\gamma$ the subgrid scale
patch field is no longer constant, and each patch grid value~$U_j$
evolves because of the coupling with its neighbours.
We construct a series solution of the \pde~(\ref{eq:lpde}) in the
coupling parameter~$\gamma$: the first order expression for the
microscopic subgrid scale field is straightforward, namely
\begin{equation}
    v_j=U_j
    +\gamma\left(\xi\mu\delta +\half\xi^2\delta^2 \right)U_j
    +c\gamma H\left( -\half r^2\xi+\rat16\xi^3 \right)\delta^2U_j
    +\Ord{\gamma^2,a^2+b^2+c^2},
    \label{eq:sub1}
\end{equation}
where the microscopic variable $\xi=(x-X_j)/H$ ranges over $|\xi|<r$\,;
accompanying these subgrid fields the grid values~$U_j$ evolve
according to standard second order discretisation (upon putting
$\gamma=1$)
\begin{equation}
    \dot U_j=
    \gamma\left( \frac1{H^2}\delta^2 U_j
    -\frac{c}{H}\mu\delta U_j \right)
    +\Ord{\gamma^2,a^2+b^2+c^2}.
    \label{eq:lowmod}
\end{equation}
See that the powers of the coupling parameter~$\gamma$ in the
\pbc~(\ref{eq:pbcs}) are chosen so that discarding terms
of~$\Ord{\gamma^p}$ results in a discrete model, such
as~(\ref{eq:lowmod}), which is of width~$2p-1$ in the grid
values~$U_j$; for example, the above model only involves $U_j$
and~$U_{j\pm 1}$\,.
Centre manifold theory~\cite[e.g.]{Carr81,Carr83b} asserts that for
small enough~$\gamma$ all neighbouring solutions are exponentially
quickly attracted to the resultant model which faithfully describes the
dynamics of the system.
Although no proof is yet available, we anticipate that the case of
interest, when $\gamma=1$\,, is small enough for this novel theoretical
support to still hold.

In the interim we demonstrate high order consistency.
We obtain models that resolve more detail of the subgrid microscopic
dynamics and its interaction with neighbouring patches by determining
higher order terms in the coupling parameter~$\gamma$.
Iteration~\cite{Roberts96a} straightforwardly generates higher order
approximations.\footnote{We use the \reduce\ computer algebra package
which has free demonstration versions available via
\url{http://reduce-algebra.com}.
Obtain our code for this problem from
\url{http://www.sci.usq.edu.au/staff/aroberts/linpbc.red}.} For
example, discarding terms~$\Ord{\gamma^3}$ the subgrid field in each
patch is modified from~(\ref{eq:sub1}) to
\begin{eqnarray}
    v_j&=&\left\{\vphantom{\frac11} 1
    +\gamma\left[\xi\mu\delta +\half\xi^2\delta^2 \right]
    +\gamma^2\left[ \rat16(\xi^3-\xi)\mu\delta^3
    +\rat1{24}(\xi^2-\xi^4)\delta^4 \right]
    \right.\nonumber\\&&\left.{}
    +c H\left[ (\gamma-\gamma^2)\rat16\xi^3\delta^2
    +\gamma^2(\rat1{60}\xi^5-\rat1{18}\xi^3)\delta^4 \right]
    +\frac bH\gamma^2\rat16\xi^3\delta^4
    \right.\nonumber\\&&\left.{}
    +r^2\left[ -cH(\gamma-\gamma^2)\rat12\xi\delta^2
    -cH\gamma^2(\rat1{12}\xi^3-\rat16\xi)\delta^4 
    -\frac bH\gamma^2\rat12\xi\delta^4 \right]
    \right.\nonumber\\&&\left.{}
    +r^4cH\rat16\xi\delta^4
    \vphantom{\frac11}\right\}U_j 
    +\Ord{\gamma^3,a^2+b^2+c^2}.
    \label{eq:sub2}
\end{eqnarray}
The first line in~(\ref{eq:sub2}) contains the leading few terms in a
universal subgrid structure for symmetric operators.  However, odd
operators, such as the advection~$c\D xu$ and the dispersion~$b\DDD
xu$, generate nontrivial subgrid structures in each patch, such as
those in the second line of~(\ref{eq:sub2}), which reflect subgrid
scale interaction of processes.  The third and fourth line of the
approximate field~(\ref{eq:sub2}) depend upon the patch size
$r=h/(2H)$\,.  But physically the subgrid scale field in each patch
should be independent of the patch size~$r$.  Although there is some
dependence in these approximations, higher orders in the coupling
parameter~$\gamma$ remove it.  For example, at the beginning of the
third line in~(\ref{eq:sub2}) see the term
$-cH(\gamma-\gamma^2)\rat12\xi\delta^2$ disappears when we set
$\gamma=1$ for the physically relevant approximation.  Similarly,
computing the next order terms in coupling parameter~$\gamma$ generates
terms, in $\gamma^3$, which cancel the $r$~dependent terms in the third
and fourth line of the subgrid field~(\ref{eq:sub2}).  Thus higher
order models push any undesirable $r$~dependence to higher orders,
thereby usefully predicting a subgrid field largely independent of the
patch size~$r$.

Simultaneously with the derivation of the subgrid field~(\ref{eq:sub2})
we determine the corresponding evolution of the macroscopic grid
values~$U_j$ for the \pde~(\ref{eq:lpde}).
Computing to higher order in the coupling parameter~$\gamma$ produces
refinements to the basic discretisation~(\ref{eq:lowmod}); for example,
here we discard terms~$\Ord{\gamma^4}$ to determine
\begin{eqnarray}
    \dot U_j
    &=&\frac1{H^2}\left( \gamma\delta^2 -\rat1{12}\gamma^2\delta^4 
    +\rat1{90}\gamma^3\delta^6 \right)U_j
    -\frac{c}{H}\left( \gamma\mu\delta 
    -\rat16\gamma^2\mu\delta^3
    +\rat1{30}\gamma^3\mu\delta^5 \right)U_j
    \nonumber\\&&{}
    -\frac{b}{H^3}\left( \gamma^2\mu\delta^3
    -\rat14\gamma^3\mu\delta^5 \right)U_j
    -\frac{a}{H^4}\left( \gamma^2\delta^4 
    -\rat16\gamma^3\delta^6 \right)U_j
    \nonumber\\&&{}
    +\Ord{\gamma^4,a^2+b^2+c^2}\,.
    \label{eq:linmod}
\end{eqnarray}
Set $\gamma=1$ to recover a model for the \pde~(\ref{eq:lpde})
supported by centre manifold theory.  Note how truncating the expansion
to different powers of coupling parameter~$\gamma$ changes the width
in~$U_j$ of the discrete model.  With the patch boundary
conditions~(\ref{eq:pbcs}) the model is independent of the patch
size~$r$.

As well as the novel dynamical systems support of exponentially quick
attractiveness and long term fidelity at finite grid size~$H$, as
mentioned earlier, another way to assess the model's relevance is to
compare the original \pde\ with the equivalent \pde\ obtained from
model~(\ref{eq:linmod}) in the limit as the macroscopic spacing $H\to
0$\,.
From~(\ref{eq:linmod}), straightforward algebra\footnote{See our
\reduce\ code from the internet.} deduces the equivalent \pde
\begin{eqnarray}&&
    \D tU =  \gamma\DD xU -\gamma c\D xU 
    -\gamma^2b\DDD xU -\gamma^2a\DDDD xU
    \nonumber\\&&{}
    +H^2\left[ (\gamma-\gamma^2)
    \left(\rat1{12}\Dn x4U -\rat16 c\Dn x3U\right)
    -(\gamma^2-\gamma^3)
    \left( \rat14b\Dn x5U +\rat16 a\Dn x6U\right) \right]
    \nonumber\\&&{}
    +H^4\left[ 
    (\rat1{360}\gamma-\rat1{72}\gamma^2+\rat1{90}\gamma^3)
    \left( \Dn x6U -3c\Dn x5U \right)
    \right.\nonumber\\&&\left.\quad{}
    -(\rat1{80}\gamma^2-\rat1{24}\gamma^3)
    \left( b\Dn x7U +2a\Dn x8U \right) \right]
    +\Ord{H^6,\gamma^4,a^2+b^2+c^2}\,.
    \label{eq:epde}
\end{eqnarray}
When the coupling parameter $\gamma=1$ the second and third lines in
the equivalent \pde~(\ref{eq:epde}) disappear and consequently the
diffusion and advection is modelled with errors of~$\Ord{H^6}$, whereas
the dispersion and the fourth-order dissipation is modelled with
errors~$\Ord{H^4}$.
Should you truncate the discretisation~(\ref{eq:linmod}) to lower
orders in coupling parameter~$\gamma$, there is less cancellation in
the equivalent \pde\ and the errors are consequently larger.
Conversely, the errors move to progressively higher orders as more
terms in the coupling parameter~$\gamma$ are retained in the centre
manifold discretisation~(\ref{eq:linmod}).
Our patch boundary conditions~(\ref{eq:pbcs}) seem to create excellent
discretisations for \pde{}s.

\section{The diffusive model is numerically stable}
\label{sec:4}

Although the \pbc{}s~(\ref{eq:pbcs}) form consistent models we need to
confirm they are numerically stable.
Indeed many other forms of \pbc{}s were tried before finding one that
was both consistent and numerically stable.
In this section we explore the gap-tooth simulations of the simple
diffusion equation
\begin{equation}
    \D tu=\DD xu\,,
    \quad\text{and $2\pi$-periodic in~$x$.}
    \label{eq:diffn}
\end{equation}
Imagine we only have access to the dynamics through a microscopic
simulator of the diffusion~(\ref{eq:diffn}), here coded by a
fine discretisation on $n$~grid points in a patch of microscopic
size~$h=rH$ and with some microscopic time step, typically $\Delta
t=10^{-6}$--$10^{-4}$.

\begin{table}
    \centering
	\caption{Growth rates~$\lambda$ of perturbations from steady state
	$u=0$\,: for diffusion~(\ref{eq:diffn}) with $m$~patches; with gap
	to patch ratio $r=0.1$\,; $n=11$ points in the microscale grid;
    and with the fourth order \pbc~(\ref{eq:pbc4}).}
    \label{tbl:spatb}
    \begin{tabular}{|r|llll|c|}
        \hline
        $m$ & \quad 1 & \quad 2,3 & \quad 4,5 & \quad 6,7 & $m+1:2m$  \\
        \hline
        4 & \E2{-12} & $ -0.946817$ & $-2.170942$ & n/a & $-99.79$  \\
        8 & \E5{-12} & $-0.996139$ & $-3.787268$ & $-7.132829$ & $-399.1$  \\
        16 & \E2{-10} & $-0.999758$ & $-3.984556$ & $-8.834269$ & $-1596.$  \\
        32 & \E{-2}{-10} & $-0.999987$ & $-3.999031$ & $-8.988851$ & $-6386.$  \\
        \hline
    \end{tabular}
\end{table}

Firstly we implement the \pbc\ that on the edge of each patch the
fine discretisation has boundary condition
\begin{equation}
    \left[\mu\delta \pm r\delta^2 
    -(\rat16-\rat12r^2)\mu\delta^3
    \mp r(\rat1{12}-\rat16r^2)\delta^4 \right]U_j
    =  H\partial_x v_j
    \quad\mbox{at }x=X_j\pm rH
    \,.
   \label{eq:pbc4}
\end{equation}
Obtain this from the first few terms of~(\ref{eq:boper}) or
equivalently from \pbc~(\ref{eq:pbcs}) by discarding $\Ord{\gamma^3}$
terms.
For the $j$th patch this \pbc\ involves macroscopic grid
values~$U_{j-2},\ldots,U_{j+2}$ only.
Then systematically perturbing each and every microscopic value from
zero, there are $mn$~such microscopic values, we numerically determined
the map of one microscopic time step.\footnote{In general, the dominant
eigenvalues of the time-stepper map may be obtained via a matrix-free
Krylov subspace iteration~\cite{Samaey03a}.
Thus for particle simulations we do not necessarily need access to all
the fine details of the microscale.} Transform the eigenvalues~$\mu$ of
this map to growth rates $\lambda=\log(\mu)/\Delta t$\,.
The $mn$~growth rates fall into $n$~groups of $m$~modes.
Each group corresponds to a microscopic internal mode of the dynamics,
roughly $\exp(\lambda_\ell t)\cos[\ell\pi(x-X_j+h/2)/h]$ for growth
rate $\lambda_\ell\approx -\ell^2\pi^2/h^2$ for $\ell=0,1,\ldots,n-1$\,.
For $\ell\geq1$ these are the rapidly decaying microscopic modes
internal to each patch seen in the initial instants of the simulations
of Figures~\ref{fig:burg3} and~\ref{fig:burg}.
The other group of $m$~modes, $\ell=0$\,, with small growth rates,
correspond to the relatively slowly evolving macroscopic modes of
interest that arise through the coupling between patches of the
microscopic dynamics.
Table~\ref{tbl:spatb} shows the leading seven growth rates, and the
magnitude of the $\ell=1$ internal growth rates, for various numbers of
patches, $m=4,8,16,32$\,.
The exact growth rates of the diffusion \pde~(\ref{eq:diffn}) are
$\lambda=-k^2$ for integer~$k$.
See in the table that as the number of patches double, the accuracy of
the growth rates of the macroscopic modes improves by a factor of
about~$16$.
This is consistent with an~$\Ord{H^4}$ method as predicted for
diffusion with \pbc~(\ref{eq:pbc4}).

\begin{table}
    \centering
	\caption{Growth rates of perturbations from steady state $u=0$ as
	for Table~\ref{tbl:spatb} but fewer points in the fine grid,
	namely $n=7$\,.}
    \label{tbl:spat7}
    \begin{tabular}{|r|llll|c|}
        \hline
        $m$ & \quad 1 & \quad 2,3 & \quad 4,5 & \quad 6,7 & $m+1:2m$  \\
        \hline
        4 & \E8{-13} & $-0.947206$ & $-2.173003$ & n/a & $-99.30$  \\
        8 & \E{-8}{-12} & $-0.996246$ & $-3.788826$ & $-7.138379$ & $-397.2$  \\
        16 & \E{-1}{-11} & $-0.999785$ & $-3.984985$ & $-8.836383$ & $-1588.$  \\
        32 & \E8{-11} & $-0.999994$ & $-3.999139$ & $-8.989397$ & $-6355.$  \\
        \hline
    \end{tabular}
\end{table}
Second, we repeat the analysis for fewer subgrid points so that the
microscopic dynamics are not resolved as well.
Table~\ref{tbl:spat7} shows the leading eigenvalues for $n=7$ points in
each patch.
There is no significant difference between Tables~\ref{tbl:spatb}
and~\ref{tbl:spat7} indicating that the microscopic resolution, the
only difference between the two, has little impact on the macroscopic
dynamics.
No growth rate is significantly positive showing the numerical method
is stable---the leading growth rate is close to zero corresponding to
conservation of material.
The other dominant growth rates rapidly approach those for diffusion.

\begin{table}
    \centering
	\caption{Growth rates of perturbations from steady state $u=0$ as
	for Table~\ref{tbl:spatb} but with the sixth order
	\pbc~(\ref{eq:pbc6}).}
    \label{tbl:spati7}
    \begin{tabular}{|r|llll|c|}
        \hline
        $m$ & \quad 1 & \quad 2,3 & \quad 4,5 & \quad 6,7 & $m+1:2m$  \\
        \hline
        4 &  \E8{-12} & $-0.982238$ & $-2.457648$ & n/a & $
        -99.79$ \\
        8 &  \E4{-11} &$-0.999677$ & $-3.928952$ & $-7.843254$ & $-399.1$ \\
        16 & \E{4}{-11} &$-1.000006$ & $-3.998708$ & $-8.967122$ & $-1596.$ \\
        32 & \E{-2}{-10} &$-1.000003$ & $-4.000023$ & $-8.999625$ & $-6386.$  \\
        \hline
    \end{tabular}
\end{table}

Lastly, consider the diffusive dynamics of~(\ref{eq:diffn}) when
connected by the \pbc\ that at $x=X_j \pm rH$
\begin{eqnarray}
    H\partial_x v_j &=&
    \left[\mu\delta \pm r\delta^2 
    -(\rat16-\rat12r^2)\mu\delta^3
    \mp r(\rat1{12}-\rat16r^2)\delta^4
    \right.\nonumber\\&&\left.\quad{}
    +(\rat1{30}-\rat18r^2+\rat1{24}r^4)\mu\delta^5
    \pm r(\rat1{90}-\rat1{36}r^2+\rat1{120}r^4)\delta^6 \right]U_j\,.
   \label{eq:pbc6}
\end{eqnarray}
Obtain this \pbc\ from the first six terms of~(\ref{eq:boper}) or
equivalently from \pbc~(\ref{eq:pbcs}) by discarding $\Ord{\gamma^4}$
terms.
Table~\ref{tbl:spati7} demonstrates that the resultant numerical scheme
is stable and has sixth order consistency for the diffusion equation.
Further, it is these \pbc{}s we used to simulate the nonlinear dynamics
of Burgers' equation~(\ref{eq:burg}) to create Figures~\ref{fig:burg3}
and~\ref{fig:burg}.

\section{Conclusion}
We achieve higher order accuracy in the gap-tooth scheme using
carefully crafted patch boundary conditions (\pbc{}s).
Analytic approximations and analysis of numerical steps in time confirm
the \pbc{}s~(\ref{eq:pbc4}) and~(\ref{eq:pbc6}) are effective.
Importantly the \pbc{}s~(\ref{eq:pbc4}) and~(\ref{eq:pbc6}) do
\emph{not} depend upon the particular \pde\ being simulated, thus the
\pbc{}s should work effectively for particle simulations for which we
do not have an algebraic microscale closure.

Further, although the predicted microscopic subgrid scale fields do
have some dependence upon the patch size~$r$, the dependence weakens,
by being pushed to higher orders in~$r$, in using higher order accuracy
patch boundary conditions.

As shown in Figures~\ref{fig:burg3} and~\ref{fig:burg}, the \pbc{}s we
recommend here appear to work well even for the nonlinear dynamics of
Burgers equation~(\ref{eq:burg}).

\bibliographystyle{plain}
\bibliography{ajr,bib,new}

\end{document}